\documentclass{amsart}

\usepackage{amsfonts}
\usepackage{amsmath}
\usepackage{amsthm}
\usepackage{amssymb}
\usepackage[english]{babel}
\usepackage{graphics}
\usepackage{latexsym}
\usepackage{longtable} 
\usepackage{mathrsfs} 
\usepackage{graphicx}
\usepackage{wrapfig} 
\usepackage[pdftex,colorlinks=true,linkcolor=blue,citecolor=magenta]{hyperref}
\usepackage{enumitem}
\usepackage{esint}

\newtheorem{theorem}{Theorem}
\newtheorem{corollary}[theorem]{Corollary}
\newtheorem{lemma}[theorem]{Lemma}
\newtheorem{proposition}[theorem]{Proposition}

\newtheorem{claim}[theorem]{Claim}
\newtheorem{example}[theorem]{Example}

\theoremstyle{definition}
\newtheorem{definition}[theorem]{Definition}
\newtheorem{remark}[theorem]{Remark}

\newcommand{\mL}{\mathcal{L}}
\newcommand{\mH}{\mathcal{H}}
\newcommand{\mF}{\mathcal{F}}

\newcommand{\mM}{\mathcal{M}}

\newcommand{\mD}{\mathcal{D}}

\newcommand{\E}{\mathrm{E}}
\newcommand{\D}{\mathrm{D}}

\newcommand{\R}{\mathbb{R}}

\newcommand{\N}{\mathbb{N}}

\newcommand{\mB}{\mathbb{B}}

\newcommand{\Q}{\mathrm{Q}}
\renewcommand{\L}{\mathrm{L}}

\newcommand{\noi}{\noindent}
\newcommand{\ms}{\medskip}

\newcommand{\al}{\alpha}

\newcommand{\ga}{\gamma}
\newcommand{\Ga}{\Gamma}
\newcommand{\de}{\delta}
\newcommand{\De}{\Delta}
\newcommand{\e}{\varepsilon}

\newcommand{\la}{\lambda}

\newcommand{\Om}{\Omega}

\newcommand{\av}{-\hspace{-10.5pt}\int}

\newcommand{\weak }{\, -\!\!\!\!-\!\!\!\!\rightharpoonup}
\newcommand{\weakstar }{ \overset{\, *_{\phantom{|}}}{{\smash{\weak }}\, } }

\newcommand{\larrow}{\longrightarrow}
\newcommand{\ot}{\otimes}

\newcommand{\LL}{\text{\LARGE$\llcorner$}}
\newcommand{\p}{\partial}
\newcommand{\sub}{\subseteq}
\newcommand{\set}{\setminus}
\newcommand{\by}{\times}

\newcommand{\ess}{\mathrm{ess}}

\renewcommand{\div}{\mathrm{div}}

\newcommand{\bt}{\begin{theorem}}\newcommand{\et}{\end{theorem}}
\newcommand{\bd}{\begin{definition}}\newcommand{\ed}{\end{definition}}
\newcommand{\bl}{\begin{lemma}}\newcommand{\el}{\end{lemma}}
\newcommand{\beq}{\begin{equation}}\newcommand{\eeq}{\end{equation}}
\newcommand{\bc}{\begin{claim}}\newcommand{\ec}{\end{claim}}
\newcommand{\bex}{\begin{example}}\newcommand{\eex}{\end{example}}
\newcommand{\bcor}{\begin{corollary}}\newcommand{\ecor}{\end{corollary}}
\newcommand{\bp}{\begin{proof}}\newcommand{\ep}{\end{proof}}

\newcommand{\BPP}{\medskip \noindent \textbf{Proof of Proposition} }

\numberwithin{equation}{section}

\begin{document}

\title[$L^\infty$ inverse source identification]{An $L^\infty$ regularisation strategy to the inverse source identification problem for elliptic equations}

\author{Nikos Katzourakis}

\address{Department of Mathematics and Statistics, University of Reading, Whiteknights, PO Box 220, Reading RG6 6AX, United Kingdom}

\email{n.katzourakis@reading.ac.uk}

  \thanks{\!\!\!\!\!\!\!\!\texttt{The author has been partially financially supported by the EPSRC grant EP/N017412/1}}
  

\date{}


\begin{abstract} In this paper we utilise new methods of Calculus of Variations in $L^\infty$ to provide a regularisation strategy to the ill-posed inverse problem of  identifying the source of a non-homogeneous linear elliptic equation, satisfying Dirichlet data on a domain. One of the advantages over the classical Tykhonov regularisation in $L^2$ is that the approximated solution of the PDE is uniformly close to the noisy measurements taken on a compact subset of the domain. 

\end{abstract}

\maketitle


\section{Introduction} \label{section1}

Let $n\in\N$ and $\Om \sub \R^n$ be a bounded domain with $C^{1,1}$ regular boundary $\p\Om$. Let also $\L$ be the linear non-divergence differential operator 
\beq \label{1.1}
\L [u]\,:=\, A :\D^2 u\,+\, b \cdot \D u \,+\, cu
\eeq
which is assumed to be uniformly elliptic with bounded continuous coefficients (and to satisfy the maximum principle):
\beq \label{1.2}
\left\{
\begin{split}
 & A \in (C^0 \cap L^\infty)(\Om;\R_s^{n\by n}) \text{ and exists } a_0>0  :\  A\!:\!\xi \ot \xi \, \geq \, a_0|\xi|^2 \text{ for }
 \\ & \text{all }\xi \in\R^n\, ;\ \ b\in (C^0\cap L^\infty)(\Om;\R^n)\, ;\ \ c \in (C^0\cap L^\infty)(\Om) \text{ and }c\leq 0 . 
\end{split}
\right.
\eeq
In the above, the notations ``$:$" and ``$\cdot$" symbolise the Euclidean inner products in the space of symmetric matrices $\R^{n\by n}_s$ and in $\R^n$ respectively, whilst $\D u =(\D_i u)_{i=1...n}$, $\D^2 u =(\D^2_{ij} u)_{i,j=1...n}$ and $\D_i\equiv \p/\p x_i$. The direct (or forward) Dirichlet problem for the above operator has the form
\beq  \label{1.3}
\left\{
\begin{array}{rl}
\L [u]\, =\, f, & \text{ in }\Om,
\\
u\, =\, g, & \text{ on }\p \Om,
\end{array}
\right.
\eeq
and asks to determine $u$, given a source $f$ and boundary data $g$. This is a classical problem which is essentially textbook material, see e.g.\ \cite[Ch.\ 9]{GT}. In particular, it is well-posed (in the sense of Hadamard) and, given $f\in L^\infty(\Om)$ and $g\in W^{2,\infty}(\Om)$, there exists a unique solution $u$ in the locally convex (Fr\'echet) space 
\beq 
\label{1.4}
\mathcal{W}^{2,\infty}_g(\Om)\,:= \bigcap_{1<p<\infty} \Big\{ u \in \big(W^{2,p}\cap W^{1,p}_g\big)(\Om) :\ \ \L[u] \in L^\infty(\Om)\Big\}.
\eeq
Note that due to the failure of the $L^p$ elliptic estimates when $p=\infty$ (see e.g.\ \cite{GM}), in general $u\not\in W^{2, \infty}(\Omega)$. Let us also note with the assumptions \eqref{1.2} on $\L$, the case of divergence operators with $C^1$ matrix coefficient $A$ is included as a special case:
\[
\L' [u]\,=\, \div(A \D u)\,+\, b \cdot \D u \,+\, cu.
\]
The {\it inverse problem} associated with \eqref{1.3} consists of the question of finding $f$, given the boundary data $g$ and some {\it partial information on the solution $u$, typically obtained through noisy (i.e.\ approximate) experimental measurements known only up to some error.} This problem is {\it severely ill-posed}, as the noisy data measured on a subset of the domain might either not be compatible with {\it any} exact solution, or even if they do, they may not suffice to determine a unique source $f$ from it.

The above inverse problem is particularly important for several applications, especially in the model case of the Laplace operator $\L=\De$ and the Poisson equation, see e.g.\ \cite{AAM, BD, EHN, I, LHY, MV, NA, SH, X, Y, YF, ZMYX}. Herein we will assume that the noisy measurements on the solution take the form
\beq
\label{1.5}
\Q[u]\, =\, q^\de \ \ \text{ on }\Gamma,
\eeq
where $\Q$ is the (nonlinear differential) {\it observation operator}
\beq
\label{1.6}
\Q[u]\,:=\, K(\cdot,u,\D u)
\eeq
with $K$ satisfying
\beq
\label{1.7}
K \in C^0(\Ga \by \R \by \R^n) \text{ and } K(x,\cdot,\cdot) \in C^1(\R\by \R^n) \text{ for any } x\in \Ga.
\eeq
Here $\Ga$ is the set on which we take measurements. It will be assumed it satisfies
\beq
\label{1.8}
\Gamma \sub \overline{\Om} \text{ is compact and }\mH^\ga(\Ga)<\infty,\text{ for some }\ga\in[0,n].
\eeq
In the above, $\mH^\ga$ denotes the Hausdorff measure of dimension $\ga$. Our general measure theory and function space notation will be either self-explanatory or otherwise standard, as e.g.\ in \cite{D,E,KV}. Finally, $q^\de \in L^\infty(\Ga,\mH^\ga)$ is the function of noisy (deterministic) measurements taken on $\Ga$, at noise level at most $\de>0$, that is
\beq 
\label{1.9}
\| q^\de - q^0 \|_{L^\infty(\Ga,\mH^\ga)} \, \leq\, \de,
\eeq
where $q^0=\Q[u^0]$ corresponds to ideal noise-free measurements of an exact solution to \eqref{1.3} with source $\L[u^0]$. 

Recapitulating, in this paper we study the following ill-posed inverse source identification problem: 
\beq
\label{1.10}
\left\{
\begin{array}{rl}
\L [u]\, =\, f\ , & \text{ in }\Om,
\\
u\, =\, g\ , & \text{ on }\p \Om,
\\
\Q[u]\, =\, q^\de, & \text{ on }\Ga.
\end{array}
\right.
\eeq
Namely, we {\it seek to specify with some selection process a suitable approximation for $f$ from measured data $q^\de$ on the compact set $\Ga$ through some observation $\Q[u]$ of the solution $u$}. Our analysis does not exclude the extreme cases $\Ga=\overline{\Om}$ (full a priori information) and $\Ga=\emptyset$ (no a priori information), although if $\Ga=\emptyset$ certain trivial modifications in the proofs are required which we do not discuss explicitly. The goal is a strategy to determine an ``optimal" best fitting solution $u^\de$ (and corresponding source $f^\de:=\L[u^\de]$) to the ill-posed problem \eqref{1.10}. In general, an exact solution may well not exist as \eqref{1.5} is a possibly incompatible pointwise constraint on $\Ga$ to the solution of \eqref{1.3} (due to the errors in measurements). On the other hand, it is not possible to have a uniquely determined source on the constraint-free region $\Om\set\Ga$, see example \ref{example}. Another popular choice in the literature for the observation operator $\Q$ consist of one of the terms in the separation of variables formula (when $\L=\De$ on rectangular domains), as e.g.\ in \cite{YF}. To the best of our knowledge, \eqref{1.10} has not been studied before in this generality.

Herein we follow an approach based on recent advances in Calculus of Variations in the space $L^\infty$ (see \cite{KM,KP1,KP2,KPa}) developed recently for functionals involving higher order derivatives. The field has been initiated in the 1960s by Gunnar Aronsson (see e.g.\ \cite{A1,A2,A3,A,AB}) and is still a very active area of research; for a review of the by-now classical theory involving scalar first order functionals we refer to \cite{K1}. To this end, we provide a {\it regularisation strategy} inspired by the classical Tykhonov regularisation strategy in $L^2$ (see e.g.\ \cite{Ki,N}), but for the following $L^\infty$ ``error" functional:
 \beq 
 \label{1.11} 
\ \ \mathrm{E}_\infty (u)\, :=\,\big \|\Q[u]-q^\de \big\|_{L^\infty(\Ga,\mH^\ga)}+\, \al \big \| \L[u] \big\|_{L^\infty( \Om)}, \ \ \ \ u\in \mathcal{W}^{2,\infty}_g(\Om),
\eeq
where $\al>0$ is a fixed regularisation parameter for the penalisation term $| \L[u] |$. In the variational language, it serves to make the functional coercive in the space. The benefit of finding a best fitting solution in $L^\infty$ is apparent: we can keep the error term $|\Q[u]-q^\de |$ due to the noise effects {\it uniformly small, not merely small on average}, which would happen if one chose to minimise the integral of the error instead of the supremum.

As it is well known to the experts of Calculus of Variations in $L^\infty$, mere (global) minimisers of supremal functionals, albeit typically easy to obtain with standard direct minimisation methods (\cite{D,FL}), they are not truly optimal and they do not share the nice ``local" minimality properties of minimisers of their integral counterparts (\cite{BN, RZ}). A popular method is to use minimisers of $L^p$ approximating functionals as $p\to \infty$ and prove appropriate convergence of such $L^p$ minimisers to a limiting $L^\infty$ minimiser. This method is fairly standard nowadays and provides a selection principle of $L^\infty$ minimisers with additional favourable properties (see e.g.\ \cite{BJW1, BP, CDP, GNP, KM, KP1}). This idea is inspired by the simple measure-theoretic fact that the $L^p$ norm (of a fixed $L^1\cap L^\infty$ function) converges to the $L^\infty$ norm of the function as $p\to \infty$.

\section{The main results}

We now give the statements of the results to be established in this paper. We will obtain {\it special} minimisers of \eqref{1.11} as limits of minimisers of
 \beq 
 \label{1.12} 
\ \ \mathrm{E}_p (u)\, :=\,\big \| |\Q[u]-q^\de|_{(p)} \big\|_{L^p(\Ga,\mH^\ga)}+\, \al \big \| |\L[u]|_{(p)}  \big\|_{L^p( \Om)}, \ \ \ \ u\in \big(W^{2,p} \cap W^{1,p}_g\big)(\Om),
\eeq
where in the above we use the normalised $L^p$ norms 
\[
\big \| f\big\|_{L^p(\Ga,\mH^\ga)}\, :=\, \left(\, {\, {\av_\Ga}} |f|^p \, \mathrm{d}\mH^\ga\right)^{1/p}, \ \ \ \big \| f\big\|_{L^p(\Om)}\, :=\, \left(\, {\, {\av_\Om}} |f|^p \, \mathrm{d}\mL^n\right)^{1/p},
\]
where the slashed integral denoting average with respect to the Hausdorff measure $\mH^\ga$ and the Lebesgue measure $\mL^n$ respectively. Further, in \eqref{1.12} $|\, \cdot \,|_{(p)}$ symbolises the following $p$-regularisation of the absolute value away from zero:
\[
|a|_{(p)}\, :=\, \sqrt{|a|^2+p^{-2}}.
\]
Let us note also that, due to our $L^p$-approximation method, as an auxiliary result we also provide an $L^p$ regularisation strategy for finite $p$ as well, which has its own merits and could be useful in itself.

\bt
[$L^\infty$ and $L^p$ regularisations of the inverse source identification problem]
\label{theorem1}
Let $\Om \sub \R^n$ be a bounded $C^{1,1}$ domain and let also $g$ be in $W^{2,\infty}(\Om)$. Suppose also the operators \eqref{1.1} and \eqref{1.6} are given, satisfying the assumptions \eqref{1.2}, \eqref{1.7}, \eqref{1.8}. Suppose further a function $q^\de \in L^\infty(\Ga,\mH^\ga)$ is given which satisfies \eqref{1.9} for $\de>0$. Let finally $\al>0$ be fixed. Then, we have the following results in relation to the problem \eqref{1.10}:

\ms

\noi {\bf (i) [Existence]} There exists a global minimiser $u_\infty\equiv u_\infty^{\al,\de} \in \mathcal{W}_{g}^{2,\infty}(\Om)$ of the functional $\E_\infty$ defined in \eqref{1.11}. In particular, we have $E_\infty(u_\infty)\leq E_\infty (v)$ for all $v\in \mathcal{W}_{g}^{2,\infty}(\Om)$ and 
\[
f_\infty \equiv f^{\al,\de}_\infty \,:= \,\L[u_\infty^{\al,\de}] \, \in L^\infty(\Om).
\]
In addition, there exist signed Radon measures 
\[
\mu_\infty \equiv \mu_\infty^{\al,\de} \, \in \, \mM(\Om), \  \  \ \nu_\infty \equiv \nu_\infty^{\al,\de} \, \in \, \mM(\Ga)
\]
such that the divergence PDE
\beq
\label{1.14}
K_r(\cdot,u_\infty,\D u_\infty)\nu_\infty \, -\, \div\big( K_p(\cdot,u_\infty,\D u_\infty)\nu_\infty  \big)\, +\, \al \L^*[\mu_\infty]\, =\, 0, \phantom{\Big|}
\eeq
is satisfied by the triplet $(u_\infty,\mu_\infty,\nu_\infty)$ in the distributional sense. In \eqref{1.14}, the operator $\L^*$ is the formal adjoint of $\L$, defined through duality, i.e.\
\[
\L^*[v]\,:=\, \div(\div (Av))\, -\, \div(bv)\,+\,cv
\]
and $K_r,K_p$ denote the partial derivatives of $K(x,r,p)$ with respect to $(r,p)\in\R\by\R^n$. Additionally, the error measure $\nu_\infty$ is supported in the closure of the subset of $\Ga$ of maximum noise, that is
\beq
\label{1.15}
\mathrm{supp} (\nu_\infty) \,\sub \, \Big\{  \big|Q[u_\infty]-q^\de\big|^\bigstar = \big\| \Q[u_\infty]-q^\de \big\|_{L^\infty(\Ga,\mH^\ga)} \Big\},
\eeq
where  $``\, (\, \cdot\, )^\bigstar \,"$ symbolises the ``essential limsup" with respect to the Radon measure $\mH^\ga\LL_\Ga$ on $\Ga$, see Proposition \ref{proposition11} that follows. If additionally the measurement function $q^\de$ is continuous on $\Ga$, \eqref{1.15} improves to
\beq
\label{1.16}
\mathrm{supp} (\nu_\infty) \,\sub \, \Big\{  \big|Q[u_\infty]-q^\de\big| = \big\| \Q[u_\infty]-q^\de \big\|_{L^\infty(\Ga,\mH^\ga)} \Big\} .
\eeq

\ms

\noi {\bf (ii) [Convergence]} For any $\al,\de>0$, the minimiser $u_\infty$ can be approximated by a family of minimisers $(u_p)_{p>n}\equiv(u_p^{\al,\de})_{p>n}$ of the respective $L^p$ functionals \eqref{1.12} and the pair of measures $(\mu_\infty,\nu_\infty) \in \mM(\Om)\by \mM(\Ga)$ can be approximated by respective absolutely continuous signed measures $(\mu_p,\nu_p)_{p>n}\equiv (\mu^{\al,\de}_p,\nu^{\al,\de}_p)_{p>n}$, as follows: 

For any $p>n$, the functional \eqref{1.12} has a global minimiser $u_p\equiv u_p^{\al,\de}$ in the space $(W^{2,p}\cap W^{1,p}_g)(\Om)$ and there exists a  sequence $p_j\larrow \infty$ as $j\to \infty$, such that
\beq
\label{1.17}
\left\{\ \ 
\begin{array}{ll}
\ \ \ \ u_p \larrow u_\infty, & \text{ in }C^{1,\kappa}(\overline{\Om}),\ \ \ \ \ \,  \text{ for any }\kappa \in (0,1),
\smallskip
\\
\D^2 u_p \weak \D^2 u_\infty, & \text{ in }L^q(\Om,\R_s^{n\by n}), \text{ for any }q \in (1,\infty), \smallskip
\end{array}
\right.
\eeq
as $p\to \infty$ along the sequence. Additionally, we have
\beq
\label{1.18}
\left\{\ \ \
\begin{split}
\nu_p\, & :=\, \frac{\big|\Q[u_p]-q^\de\big|^{p-2}_{(p)} \big(\Q[u_p]-q^\de\big) }{ \mH^\ga(\Ga) \, \big \| |\Q[u_p]-q^\de|_{(p)} \big\|^{p-1}_{L^p(\Ga,\mH^\ga)}}\mH^\ga \LL_\Ga \, \weakstar \, \nu_\infty,  \ \ \ \text{ in }\mM(\Ga),
\\
\mu_p\, & :=\, \frac{| \L[u_p] |^{p-2}_{(p)} \,  \L[u_p] }{ \mL^n(\Om) \, \big \| | \L[u_p] |_{(p)} \big\|^{p-1}_{L^p(\Om)}}\mL^n \LL_\Om \, \weakstar \, \mu_\infty,  \hspace{48pt} \text{ in }\mM(\Om),
\end{split}
\right.
\eeq
as $p\to \infty$ along the sequence. Further, for each $p>n$,  the triplet $(u_p,\mu_p,\nu_p)$ solves the equation
\beq
\label{1.19}
K_r(\cdot,u_p,\D u_p)\nu_p \, -\, \div\big( K_p(\cdot,u_p,\D u_p)\nu_p  \big)\, +\, \al \L^*[\mu_p]\, =\, 0,  \phantom{\Big|}
\eeq
in the distributional sense.

\ms

\noi {\bf (iii) [$L^\infty$ error estimates]} For any exact solution $u^0 \in \mathcal{W}^{2,\infty}_g(\Om)$ of \eqref{1.10} (with $f=\L[u^0]$ and $\Q[u^0]=q^0$) corresponding to measurements with zero noise, we have the estimate: 
\beq
\label{1.20}
\Big\| \Q[u^{\al,\de}_\infty] - \Q[u^0] \Big\|_{L^\infty(\Ga,\mH^\ga)}  \leq\, 2\de + \al \, \| \L[u^0]\|_{L^\infty(\Om)},
\eeq
for any $\al,\de>0$.

\ms

\noi {\bf (iv) [$L^p$ error estimates]}  For any exact solution $u^0 \in (W^{2,p}\cap W^{1,p}_g)(\Om)$ of \eqref{1.10} (with $f=\L[u^0]$ and $\Q[u^0]=q^0$) corresponding to measurements with zero noise and for $p>n$, we have the estimate: 
\beq
\label{1.24}
\Big\| \Q[u^{\al,\de}_p] - \Q[u^0] \Big\|_{L^p(\Ga,\mH^\ga)}  \leq\, 2\de + \al \, \| \L[u^0]\|_{L^p(\Om)},
\eeq
for any $\al,\de>0$.
\et 
The estimate in part {\bf (iv)} above is useful if we have merely that $\L[u^0] \in L^p (\Om)$ for $p<\infty$ (namely when perhaps $\L[u^0] \not\in L^\infty(\Om)$).

\ms

\section{Discussion and auxiliary results}

We begin by providing some clarifications regarding Theorem \ref{theorem1}.

\begin{remark} \label{remark4} {\bf (i)} We note that in \eqref{1.14} the distributional meaning of this PDE is
\[
\int_\Ga K_r(\cdot,u_\infty,\D u_\infty) \,\phi\, \mathrm{d} \nu_\infty \, +\, \int_\Ga K_p(\cdot,u_\infty,\D u_\infty)\cdot \D \phi \, \mathrm{d} \nu_\infty \, +\, \al \int_\Om \L[\phi] \, \mathrm{d} \mu_\infty\, =\, 0,
\]
for all test functions $\phi \in C^2_c(\Om)$. Therefore, in fact the equation \eqref{1.14} is valid in the smaller space of {\it second order distributions}:
\[
\mD^{-2}(\Om)\,:=\, \big(C^2_c(\Om)\big)^*.
\]
Additionally, since the measure $\nu_\infty$ is supported in the compact set $\Ga$, by extending $\nu_\infty$ on $\Om \set \Ga$ by zero (i.e.\ by identifying $\nu_\infty$ with the restriction $\nu_\infty \LL _\Ga$), we may rewrite \eqref{1.14} as
\[
\int_\Om \Big( K_r(\cdot,u_\infty,\D u_\infty) \,\phi \, +\, K_p(\cdot,u_\infty,\D u_\infty)\cdot \D \phi \Big)\, \mathrm{d} \nu_\infty \, +\, \al \int_\Om \L[\phi] \, \mathrm{d} \mu_\infty\, =\, 0,
\]
for all $\phi \in C^2_c(\Om)$. 

\ms

\noi {\bf (ii)} In index form, the definition of the formal adjoint can be written as
\[
\L^*[v]\,=\, \sum_{i,j=1}^n\D^2_{ij}(A_{ij}v)\, -\, \sum_{k=1}^n\D_k(b_k\, v)\,+\,cv
\]
and the distributional interpretation of $\L^*$ through duality is
\[
\langle \L^*[v],\phi\rangle\,=\, \int_\Om \bigg(\sum_{i,j=1}^n(\D^2_{ij}\phi)(A_{ij}v)\, +\, \sum_{k=1}^n(\D_k\phi)( b_k v)\,+\,\phi cv \bigg)\, \mathrm{d} \mL^n,
\]
for all $\phi \in C^2_c(\Om)$. In a similar vein, the distributional interpretation of \eqref{1.19} is
\[
\int_\Ga  K_r(\cdot,u_p,\D u_p) \,\phi \, \mathrm{d} \nu_p \, +\, \int_\Ga K_p(\cdot,u_p,\D u_p)\cdot \D \phi \, \mathrm{d} \nu_p \, +\, \al \int_\Om \L[\phi] \, \mathrm{d} \mu_p\, =\, 0,
\]
for all $\phi \in C^2_c(\Om)$. By taking into account that the measures $\mu_p,\nu_p$ as given by \eqref{1.18} are in fact absolutely continuous with respect to the Lebesgue and the Hausdorff measure respectively, the above is in fact equivalent to
\[
\begin{split}
\, {\av_\Ga} \Big( K_r(\cdot,u_p,\D u_p) \,\phi \, +\, K_p(\cdot,u_p,\D u_p)\cdot \D \phi \Big) \frac{\big|\Q[u_p]-q^\de\big|^{p-2}_{(p)} \big(\Q[u_p]-q^\de\big)}{\big \| |\Q[u_p]-q^\de|_{(p)} \big\|^{p-1}_{L^p(\Ga,\mH^\ga)}} \, \mathrm{d} \mH^\ga \, 
\\
+\ \al \, {\, {\av_\Om}} \L[\phi] \frac{| \L[u_p] |^{p-2}_{(p)} \,  \L[u_p]  }{ \big \| | \L[u_p] |_{(p)} \big\|^{p-1}_{L^p(\Om)}}  \, \mathrm{d} \mL^n\, =\, 0,
\end{split}
\]
for all $\phi \in C^2_c(\Om)$. 

\ms

\noi {\bf (iii)} Since we only prescribe boundary conditions $u=g$ on $\p\Om$ but impose no condition on the gradient (as opposed to e.g.\ \cite{KM}, wherein an $L^\infty$ minimisation problem was considered by imposing $\D u=\D g$ on $\p\Om$ additionally to $u=g$ on $\p\Om$), we therefore have ``natural boundary conditions" for the gradient on $\p\Om$. We will make no particular further use of this observation.
\end{remark}

The following two results are consequences of our main theorem.

\begin{corollary}[Rates of convergence] \label{Cor2}  In the setting of Theorem \ref{theorem1}, in the case that $\Q[u]:=u$, the estimates \eqref{1.20}-\eqref{1.24} for the $L^\infty$ and the $L^p$ minimisers can be improved to the linear rates of convergence
\beq
\label{1.21}
\big\| u^{\al,\de}_\infty - u^0 \big\|_{L^\infty(\Ga,\mH^\ga)}  \leq\, 2\de + \al \, \| \L[u^0]\|_{L^\infty(\Om)} \ \ \text{ as }\al,\de\to0,
\eeq
if $\L[u^0] \in L^\infty(\Om)$, and 
\beq
\label{1.25}
\big\| u^{\al,\de}_p - u^0 \big\|_{L^p(\Ga,\mH^\ga)}  \leq\, 2\de + \al \, \| \L[u^0]\|_{L^p(\Om)}  \ \ \text{ as }\al,\de\to0,
\eeq
if $\L[u^0] \in L^p(\Om)$ for $p<\infty$. 
\end{corollary}

\begin{corollary} \label{Cor3} In the setting of Theorem \ref{theorem1}, we have 
\[
L^*[\mu_p] \, =\, 0 \ \  \text{ in } \Om\set\Ga
\]
in the distributional sense, for any $p\in(1,\infty]$. In particular, for $p<\infty$ we have
\[
\L^*\Big(\big| \L[u_p] \big|^{p-2}_{(p)} \, \L[u_p]\Big) \, =\, 0  \ \  \text{ in }\Om\set \Ga,
\]
in the distributional sense.
\end{corollary}

Corollary \ref{Cor3} expresses the fact that on the subset where we have no a priori information on the solution generating the source (and hence no constraint on the PDE), then one can select a solution whose's source is associated with a solution of the dual homogeneous problem $\L^*[\mu_\infty]=0$.

\begin{remark} Possible choices for the observation operator $\Q$ which are popular in the literature, are the following:

\begin{itemize}
\item $\Q[u]:=u(x,c)$, for $n=2$ and $\Om=(a,b)\by (c,d)$ being a rectangular domain (i.e., one of the products in the separation of variables when $\L=\De$). This implies that \eqref{1.20} simplifies to
\[
\big\| u^{\al,\de}_\infty(\cdot,c) - u^0(\cdot,c) \big\|_{L^\infty((a,b),\mH^1)}  \leq\, 2\de + \al \, \| \L[u^0]\|_{L^\infty((a,b)\by (c,d))} \ \ \text{ as }\al,\de\to0,
\]
and similarly for its $L^p$-counterpart.
\ms

\item $\Q[u]:=\D u \cdot n$, where $n$ is the outer normal vector on $\p\Om$. In this case, \eqref{1.21} simplifies to
\[
\big\| n\cdot \big(\D u^{\al,\de}_\infty - \D u^0 \big)\big\|_{L^\infty(\p\Om,\mH^{n-1})}  \leq\, 2\de + \al \, \| \L[u^0]\|_{L^\infty(\Om)} \ \ \text{ as }\al,\de\to0,
\]
and similarly for its $L^p$-counterpart.
\end{itemize}
We remark that, due to the ill-posed nature of the problem, in general it is not possible to obtain an estimate on $\Om\set \Ga$, see Example \ref{example} that follows.
\end{remark}

\begin{remark}[On the source determination] We would like to point out explicitly that our result allows to construct the next putative $(\al,\de)$-dependant source for the inverse problem associated with \eqref{1.10}:
\[
f_\infty^{\al,\de} \,:= \,\L[u_\infty^{\al,\de}] \, \in L^\infty(\Om),
\]
where $u_\infty^{\al,\de}$ is the $\E_\infty$-minimiser of the regularised error in Theorem \ref{theorem1}. The natural question then arises regarding when this approximate source converges to the actual source as $\al,\de \to 0$. Unfortunately, one can not generally improve parts (iii) and (iv) of the theorem to hold on $\Om\set \Ga$ and we are bound to have convergence ``through $\Q$" on $\Ga$ only. {\it The main obstruction is that, in general, one cannot determine a unique source from the data, unless the set $\Ga$ is considerably large and the operator $\Q$ is relatively special.} 
\end{remark}

We now give an example showing that, in general, it is not possible to have a uniquely determined source on the constraint-free region $\Om\set\Ga$. In fact, if $\Ga \sub \p\Om$ (in which case $\Om\set \Ga=\emptyset$), then there is complete indeterminacy of the (solution and the) source. This is not an issue of regularity of neither the solution nor the source, as in the example below all admissible sources are equally smooth because they are perturbations parametrised by harmonic functions.

\begin{example}[See also \cite{BD}] \label{example} Let us choose
\[
\L=\De, \ \ \Q[u]= \mathrm n\cdot \D u= \frac{\p u}{\p \mathrm n}, \ \ \mathrm n\text{ the outer normal vector on }\p\Om, \ \ \Ga=\p\Om.
\]
Then, we have $\Om\set \Ga=\emptyset$ and the Dirichlet problem \eqref{1.10} becomes
\[
\left\{
\begin{array}{rl}
\De u \, =\, f\ , & \text{ in }\Om, \ms
\\
u\, =\, g\ , & \text{ on }\p \Om,
\\
\dfrac{\p u}{\p  \mathrm n} \, =\, q^\de, & \text{ on }\p \Om.
\end{array}
\right.
\]
Let $h$ be {\bf any} harmonic function on $\Om$. Let also $w$ be the unique solution to
\[
\left\{
\begin{array}{rl}
\De w \, =\, h\ , & \text{ in }\Om, \ms
\\
w\, =\, 0\ , & \text{ on }\p \Om,
\end{array}
\right.
\]
and let $v$ be the unique solution to
\[
\left\{
\begin{array}{rll}
\De^2 v \! \!  \!  &= \, f\ , & \text{ in }\Om, \ms
\\
v \! \!  \!   &= \, g\ , & \text{ on }\p \Om,
\\
\dfrac{\p v}{\p  \mathrm n} \! \!  \!   &= \, q^\de-\dfrac{\p w}{\p  \mathrm n}, & \text{ on }\p \Om.
\end{array}
\right.
\]
It follows that any source of the form $f:=h+\De v$ is associated with a solution $u:=v+w$ to the Dirichlet problem. Indeed, we have
\[
\begin{split}
u|_{\p\Om} &=\, v|_{\p\Om} +\, w|_{\p\Om}=\, v|_{\p\Om} +\, 0\, =\, g,
\\
\dfrac{\p u}{\p  \mathrm n} \Big|_{\p\Om} & =\, \dfrac{\p v}{\p  \mathrm n} \Big|_{\p\Om} +\, \dfrac{\p w}{\p  \mathrm n} \Big|_{\p\Om}=\, \bigg(q^\de-\dfrac{\p w}{\p  \mathrm n} \bigg) + \, \dfrac{\p w}{\p  \mathrm n} =\, q^\de,
\\
\De u \, &=\, \De v \, + \De w \, =\,  \De v \, + h \, =\, f. 
\end{split}
\]
This happens because the boundary data $u=g$  on $\p\Om$ and $\p u /\p n =q^\de$ on $\p\Om$ can only determine a unique biharmonic function $v$ in $\Om$ with $\De^2v=0$. 
\end{example}

Given that, as the above example certifies, one cannot determine a unique source on $\Om\set\Ga$, the result that follows provides some sufficient conditions regarding when the approximate source $f_\infty^{\al,\de}$ converges to the actual source of the problem at least on the set $\Ga$ as $\al,\de \to 0$.

\begin{corollary}[Approximation of the actual source] In the setting of Theorem \ref{theorem1} and Corollary \ref{Cor2}, let $f_\infty^{\al,\de}=\L[u_\infty^{\al,\de}]$ be the approximate source and let $\L[u^0]$ be the actual source. Suppose that the coefficients $A,b$ of $\L$ are $W^{2,\infty},W^{1,\infty}$ respectively. Suppose further that $\Q[u]=u$, $\Ga$ has non-empty interior $\Ga^\circ$ and that $\ga=n$. Then, we have that
\[
f_\infty^{\al,\de} \weak \L[u^0] \ \  \text{ in } \mD'(\Ga^\circ),
\]
distributionally as $\al,\de \to 0$. In fact, we have the following estimate which implies strong convergence in the dual Sobolev space $W^{-2,\infty}(\Ga^\circ)=(W^{2,1}_0(\Ga^\circ))^*$:
\[
\begin{split}
 \big\| f_\infty^{\al,\de}-\L[u^0]\big\|_{W^{-2,\infty}({\Ga^\circ})}  \, \leq \, C\Big(2 \de + \al \| \L[u^0]\|_{L^\infty(\Om)}\Big) ,
\end{split}
\]
where $C$ depends only on the coefficients $A,b,c$ of $\L$.
\end{corollary}

\bp By Corollary \ref{Cor2} and our assumptions on $\Q,\Ga,\ga$, we have 
\[
\big\| u^{\al,\de}_\infty - u^0 \big\|_{L^\infty(\Ga^\circ)}  \leq\, 2\de + \al \, \| \L[u^0]\|_{L^\infty(\Om)}.
\]
Fix $\phi \in C^2_c(\Ga^\circ)$. Then, by our assumption on $\L$ we have
\[
\begin{split}
\bigg| \int_{\Ga^\circ} \big( f_\infty^{\al,\de}-\L[u^0]\big)\,\phi \, \mathrm d \mL^n \bigg| \, &=\, \bigg| \int_{\Ga^\circ} \L\big[ u_\infty^{\al,\de}-u^0\big]\,\phi \, \mathrm d \mL^n \bigg| 
\\
& =\,\bigg|  \int_{\Ga^\circ} \big( u_\infty^{\al,\de}-u^0\big)\,(\L^*[\phi]) \, \mathrm d \mL^n \bigg| 
\\
& \leq \, \big\| u_\infty^{\al,\de}-u^0\big\|_{L^\infty(\Ga^\circ)} \big\| \L^*[\phi] \big\|_{L^1(\Ga^\circ)} 
\\
& \leq \, \Big(2\de + \al \, \| \L[u^0]\|_{L^\infty(\Om)}\Big) \big\| \L^*[\phi] \big\|_{L^1(\Ga^\circ)}
\\
& \leq \, C\Big(2 \de + \al \| \L[u^0]\|_{L^\infty(\Om)}\Big) \big\| \phi \big\|_{W^{2,1}_0(\Ga^\circ)},
\\
\end{split}
\]
where $C$ depends only on the coefficients of $\L$. Hence, $f_\infty^{\al,\de} \weak \L[u^0]$  in $\mD'(\Ga^\circ)$ as $\al,\de\to 0$, as claimed. The stronger convergence $f_\infty^{\al,\de} \larrow \L[u^0]$ in the Sobolev space is a consequence of the definition of the operator norm on $W^{-2,\infty}(\Ga^\circ)$.
\ep

The following result studies the ``concentration measures" of the approximate $L^p$ minimisation problems as $p\to \infty$. Note that we are actually using ``$k$" instead of  ``$p$" to avoid confusion, as we will later apply it to a certain subsequence $(p_k)_1^\infty$.

\begin{proposition}[The essential limsup] \label{proposition11} Let $X \sub \R^n$ be a Borel set, endowed with the induced Euclidean topology and let also $\nu \in \mM(X)$ be a positive finite Radon measure on $X$. For any $f\in L^\infty(X,\nu)$, we define the function $f^\bigstar \in L^\infty(X,\nu)$ by setting
\[
f^\bigstar(x) \, := \, \lim_{\e\to 0} \bigg(\nu-\underset{y \in \mB_\e(x)}{\ess\,\sup} \, f(y)\bigg)
\]
and we call $f^\star$ {\bf the $\nu$-essential limsup of $f$}. In the above, $\mB_\e(x)$ symbolises the open ball of radius $\e$ centred at $x\in X$ with respect to the induced topology. Then, we have:

\ms 

\noi {\bf (i)} It holds that $f \leq f^\bigstar$, $\nu$-a.e.\ on $X$.

\ms

\noi {\bf (ii)} It holds that $f^\bigstar$ is upper semicontinuous on $X$, namely
\[
\limsup_{X \ni y\to x} f^\bigstar(y) \, \leq \, f^\bigstar(x), \ \ \ x\in X. 
\]
{\bf (iii)} $f^\bigstar$ gives a pointwise meaning to the essential supremum on $X$, in the sense
\[
\sup_X f^\bigstar \, =\, \nu-\underset{X}{\ess\,\sup} \, f.
\]
\end{proposition}

The following result studies what we call ``concentration measures" of the approximate $L^p$ minimisation problems as $p\to \infty$. Note that we are using dumb variable ``$k$" instead of  ``$p$" to avoid confusion, as we will later apply it to a certain subsequence $(p_k)_1^\infty$.

\begin{proposition}[$L^k$ concentration measures as $k\to\infty$] \label{proposition12} Let $X$ be a compact metric space, endowed with a non-negative finite Borel measure $\nu$ which attaches positive values to any non-empty open set on $X$. Consider a sequence $(f_k)_1^\infty \sub L^\infty(X,\nu)$ and consider the sequence of absolutely continuous signed Radon measures $(\nu_k)_1^\infty \sub \mM(X)$, given by:
\[
\nu_k\,:=\, \frac{1}{ \nu(X)} \frac{\big(|f_k|_{(k)}\big)^{k-2} f_k }{ \big\| |f_k|_{(k)} \big\|^{k-1}_{L^k(X,\nu)}}\, \nu, \ \ \ k\in \N,
\]
where $|\cdot|_{(k)}=(|\cdot|^2+k^{-2})^{1/2}$. Then:
\ms

\noi {\bf (i)} There exists a subsequence $(k_i)_1^\infty$ and a limit measure $\nu_\infty \in \mM(X)$ such that
\[
\nu_k \, \weakstar \, \nu_\infty \ \ \text{ in }\mM(X),
\]
as $k_i\to \infty$.

\ms

\noi {\bf (ii)} If there exists $f_\infty \in L^\infty(X,\nu)\set \{0\}$ such that
\[
\sup_X |f_k -f_\infty| \larrow 0\ \ \text{ as }k\to\infty,
\]
then the limit measure is supported in the set where (the $\nu$-essential limsup of) $|f_\infty|$ equals $\|f_\infty\|_{L^\infty(X,\nu)}$:
\[
\mathrm{supp}(\nu_\infty) \, \sub \, \Big\{ |f_\infty|^\bigstar= \|f_\infty\|_{L^\infty(X,\nu)} \Big\}.
\]

\noi {\bf (iii)} If additionally to the assumptions of {\bf(ii)} the modulus $|f_\infty|$ of the uniform limit $f_\infty$ is continuous on $X$, then the following stronger assertion holds true:
\[
\mathrm{supp}(\nu_\infty) \, \sub\, \Big\{ |f_\infty|= \|f_\infty\|_{L^\infty(X,\nu)} \Big\}.
\]
\end{proposition}

\ms

\section{Proofs} \label{section2}

Herein we establish Theorem \ref{theorem1} and its corollaries, together with the auxiliary results Propositions \ref{proposition11}-\ref{proposition12}. The proof of Theorem \ref{theorem1} consists of several lemmas. We note that some of the details might be standard to the experts of Calculus of Variations, but we do provide most of the niceties for the sake of completeness and for the convenience of the reader.

\begin{lemma} \label{lemma1} For any $p>n$ and fixed $\al,\de>0$, the functional \eqref{1.12} has a (global) minimiser $u_p \in (W^{2,p}\cap W^{1,p}_g)(\Om)$:
\[
\E_p(u_p) \, =\, \inf\Big\{\E_p(v)\ : \ v\in (W^{2,p}\cap W^{1,p}_g)(\Om) \Big\}.
\]
\end{lemma}

\bp Since $g\in W^{2,\infty}(\Om)$ (and in particular because $g,\D g$ are continuous on $\Ga$ and therefore $\mH^\ga$-measurable by identification with their precise Lebesgue representatives reconstructed through limits of average values), by the H\"older inequality and our assumption we have the a priori bound
\[
\begin{split}
\E_p(g)\, \leq \,&\, \E_\infty(g)\\
 \leq \, & \, \|q^\de\|_{L^\infty(\Ga,\mH^\ga)}\, +\  \|K(\cdot,g,\D g) \|_{L^\infty(\Ga,\mH^\ga)}
\\
& +\al \Big(\|A\|_{L^\infty(\Om)}+\|b\|_{L^\infty(\Om)}+\|c\|_{L^\infty(\Om)} \Big) \|g\|_{W^{2,\infty}(\Om)}
\\
<\, &\, \infty.
\end{split}
\]
Hence, 
\[
0\, \leq\, \inf\Big\{\E_p(v)\ : \ v\in (W^{2,p}\cap W^{1,p}_g)(\Om) \Big\}\, \leq\, \E_\infty(g) \,<\,\infty.
\]
Further, $\E_p$ is coercive in the space $(W^{2,p}\cap W^{1,p}_g)(\Om)$: indeed, by the $L^p$ elliptic estimates for linear second order equations with measurable coefficients \cite[Ch.\ 9]{GT}, by our assumptions on $\L$ and the H\"older inequality we have
\[
\begin{split}
\E_p(v) \, &\geq\, \al \| \L[v] \|_{L^p(\Om)}
\\ 
& \geq \, \frac{\al}{C(p,A,b,c)}\Big( \|v\|_{W^{2,p}(\Om)}  -  \|g\|_{W^{2,p}(\Om)} \Big)
\\
& \geq \, \frac{\al}{C(p,A,b,c)}\Big( \|v\|_{W^{2,p}(\Om)}  -  \|g\|_{W^{2,\infty}(\Om)} 
\Big)
\end{split}
\]
for some $C=C(p,A,b,c)>0$ and any $v\in (W^{2,p}\cap W^{1,p}_g)(\Om)$. Let $(u_p^m)_1^\infty$ be a minimising sequence of $\E_p$:
\[
\E_p(u_p^m) \larrow \inf\Big\{\E_p(v)\ : \ v\in (W^{2,p}\cap W^{1,p}_g)(\Om) \Big\},
\]
as $m\to \infty$. Then, by the above estimates, we have the uniform bound
\[
 \|u_p^m\|_{W^{2,p}(\Om)}\, \leq\, C
\]
for some $C>0$ depending on $p$ but independent of $m\in\N$. By standard weak and strong compactness arguments in Sobolev spaces, there exists a subsequence $(u_p^{m_k})_1^\infty$ and a function $u_p\in (W^{2,p}\cap W^{1,p}_g)(\Om)$ such that, along this subsequence we have
\[
\left\{\ \ 
\begin{array}{ll}
u_p^m \larrow u_p, & \text{ in }L^{p}(\Om),
\smallskip
\\
\D u_p^m \larrow \D u_p, & \text{ in }L^{p}(\Om,\R^n),
\smallskip
\\
\D^2 u_p^m \weak \D^2 u_p, & \text{ in }L^p(\Om,\R_s^{n\by n}), 
\end{array}
\right.
\]
as $m_k\to \infty$. Additionally, since $p>n$, by the regularity of the boundary we have the compact embedding $W^{2,p}(\Om)\Subset C^{1,k}(\overline{\Om})$ as a consequence of the Morrey estimate. Hence,
\[
u_p^m \larrow u_p \ \ \text{ in }C^{1,\kappa}(\overline{\Om}), \text{ for }\kappa \in \left(0,1-\frac{n}{p}\right),
\]
as $m_k\to \infty$. The above modes of convergence and the continuity of the function $K$ defining the operator $\Q$ imply that $\Q[u^m_p] \larrow \Q[u_p]$ uniformly on $\Ga$ as  $m_k \to \infty$. Therefore,
\[
\big\||\Q[u^m_p]-q^\de|_{(p)} \big\|_{L^p(\Ga,\mH^\ga)} \larrow \big\||\Q[u_p]-q^\de|_{(p)} \big\|_{L^p(\Ga,\mH^\ga)}
\]
as $m_k \to \infty$. Additionally, by the linearity of the operator $\L$ and because its coefficients are $L^\infty$, we have that
\[
\L[u^m_p] \weak \, \L[u_p]\ \ \text{ in }L^{p}(\Om),
\]
as $m_k \to \infty$. Since the functional 
\[
\big\|| \cdot |_{(p)} \big\|_{L^p(\Om)} \ : \ \ L^{p}(\Om) \larrow \R
\] 
is convex on this reflexive space and also it is strongly continuous, it is weakly lower semi-continuous and therefore
\[
\big\|| \L[u_p] |_{(p)} \big\|_{L^p(\Om)}\, \leq \, \liminf_{k\to\infty} \, \big\|| \L[u^{m_k}_p]|_{(p)} \big\|_{L^p(\Om)}.
\]
By putting all the above together, we see that
\[
\E_p(u_p) \, \leq \, \liminf_{k\to\infty} \, \E_p(u^{m_k}_p) \,\leq\, \inf\Big\{\E_p(v)\ : \ v\in (W^{2,p}\cap W^{1,p}_g)(\Om) \Big\},
\]
which concludes the proof.
\ep

\begin{lemma} \label{lemma2} For any $\al,\de>0$, there exists a (global) minimiser $u_\infty \in \mathcal{W}^{2,\infty}_g(\Om)$ and a sequence of minimisers $(u_{p_i})_1^\infty$ of the respective $\E_p$-functionals constructed in Lemma \ref{lemma1}, such that \eqref{1.17} holds true.
\end{lemma}

\bp For each $p>n$, let $u_p \in (W^{2,p}\cap W^{1,p}_g)(\Om)$ be the minimiser of $\E_p$ given by Lemma \ref{lemma1}. For any fixed $q \in (n,\infty)$ and $p\geq q$, the H\"older inequality and the minimality property imply the estimates
\[
\E_q(u_p)\, \leq\, \E_p(u_p)\, \leq\, \E_p(g)\, \leq\, \E_\infty(g)\, <\, \infty.
\]
By the coercivity of $\E_q$ in the space $(W^{2,q}\cap W^{1,q}_g)(\Om)$, we have the estimate
\[
\begin{split}
\E_q(u_p) \, \geq\, \frac{\al}{C(q,A,b,c)}\Big( \| u_p \|_{W^{2,q}(\Om)}  -  \|g\|_{W^{2,\infty}(\Om)}
\Big) ,
\end{split}
\]
which implies 
\[
\sup_{p\geq q} \, \| u_p \|_{W^{2,q}(\Om)}\, \leq\, C
\]
for some $C>0$ depending on $q$, the coefficient of $\L$ and $\al$. By a standard diagonal argument, for any sequence $(p_i)_1^\infty$ with $p_i \larrow \infty$ as $i\to \infty$, there exists a function
\[
u_\infty \in \bigcap_{n<q<\infty}(W^{2,q}\cap W^{1,q}_g)(\Om)
\]
and a subsequence (denoted again by $(p_i)_1^\infty$) along which \eqref{1.17} holds true. It remains to show that $\L[u_\infty] \in L^\infty(\Om)$ (which would guarantee membership in the space $\mathcal{W}^{2,\infty}_g(\Om)$) and that $u_\infty$ is in fact a minimiser of $\E_\infty$ over the same space. To this end, note that for any fixed $q\in (n,\infty)$ and $p\geq q$, we have
\[
\E_q(u_p)\, \leq\, \E_p(u_p)\, \leq\, \E_p(v)\, \leq\, \E_\infty(v)
\]
for any $v\in \mathcal{W}^{2,\infty}_g(\Om)$. By the weak lower semi-continuity of $\E_q$ in the space $(W^{2,q}\cap W^{1,q}_g)(\Om)$ demonstrated in Lemma \ref{lemma1}, we have
\[
\E_q(u_\infty)\, \leq\, \liminf_{i\to \infty} \, \E_q(u_{p_i})\, \leq\, \E_\infty(v),
\]
for any $v\in \mathcal{W}^{2,\infty}_g(\Om)$. The particular choice $v:=g$ in the above estimate gives the bound
\[
\al \| \L[u_\infty]\|_{L^q(\Om)}\, \leq\, \E_q(u_\infty) \, \leq\, \E_\infty(g).
\]
By letting $q\to \infty$ in the last two estimates above, we obtain that $\L[u_\infty] \in L^\infty(\Om)$ and that
\[
\E_\infty(u_\infty)\,\, \leq\, \inf\Big\{ \E_\infty(v) \ : \ v \in \mathcal{W}^{2,\infty}_g(\Om)\Big\},
\]
as desired.
\ep

\begin{lemma} \label{lemma3} For any $\al,\de>0$ and $p>n$, consider the minimiser $u_p \in (W^{2,p}\cap W^{1,p}_g)(\Om)$ of the functional $\E_p$ constructed in Lemma \ref{lemma1}. Consider also the signed Radon measures $\mu_p \in \mM(\Om)$ and $\nu_p \in \mM(\Ga)$, defined as in \eqref{1.18}: 
\[
\begin{split}
\nu_p\, & :=\, \frac{\big|\Q[u_p]-q^\de\big|^{p-2}_{(p)} \big(\Q[u_p]-q^\de\big) }{ \mH^\ga(\Ga) \, \big \| |\Q[u_p]-q^\de|_{(p)} \big\|^{p-1}_{L^p(\Ga,\mH^\ga)}} \, \mH^\ga \LL_\Ga,
\\
\mu_p\, & :=\, \frac{| \L[u_p] |^{p-2}_{(p)} \,  \L[u_p]  }{ \mL^n(\Om) \, \big \| | \L[u_p] |_{(p)} \big\|^{p-1}_{L^p(\Om)}} \, \mL^n \LL_\Om .
\end{split}
\]
Then, the triplet $(u_p,\mu_p,\nu_p)$ satisfies the PDE \eqref{1.19} in the distributional sense. In fact, the following stronger assertion holds: we have
\[
\begin{split}
\, {\av_\Ga} \Big( K_r(\cdot,u_p,\D u_p) \,\phi \, +\, K_p(\cdot,u_p,\D u_p)\cdot \D \phi \Big) \frac{\big|\Q[u_p]-q^\de\big|^{p-2}_{(p)} \big(\Q[u_p]-q^\de\big)}{\big \| |\Q[u_p]-q^\de|_{(p)} \big\|^{p-1}_{L^p(\Ga,\mH^\ga)}} \, \mathrm{d} \mH^\ga \, 
\\
+\ \al \, {\, {\av_\Om}} \L[\phi] \frac{| \L[u_p] |^{p-2}_{(p)} \,  \L[u_p]  }{ \big \| | \L[u_p] |_{(p)} \big\|^{p-1}_{L^p(\Om)}}  \, \mathrm{d} \mL^n\, =\, 0,
\end{split}
\]
for all $\phi \in W^{2,p}_0(\Om)$. 
\end{lemma}

\bp We involve a standard Gateaux differentiability argument. Let us begin by checking that $\mu_p,\nu_p$ indeed define measures when $u_p \in W^{2,p}(\Om)$. Indeed, by the H\"older inequality, we have the total variation estimates
\[
\begin{split}
\|\nu_p\|(\Ga)\, & \leq \, \Big(\big \| |\Q[u_p]-q^\de|_{(p)} \big\|_{L^p(\Ga,\mH^\ga)}\Big)^{1-p} \, {\av_\Ga} \big|\Q[u_p]-q^\de\big|^{p-1}_{(p)}  \, \mathrm{d} \mH^\ga
\\
&\leq\, \Big(\big \| |\Q[u_p]-q^\de|_{(p)} \big\|_{L^p(\Ga,\mH^\ga)}\Big)^{1-p} \left(\, \, {\av_\Ga} \big|\Q[u_p]-q^\de\big|^{p}_{(p)}  \, \mathrm{d} \mH^\ga \right)^{\!\!\frac{p-1}{p}}
\\
& = \, 1
\end{split}
\]
and similarly
\[
\begin{split}
\|\mu_p\|(\Om)\, & \leq \, \Big(\big \| | \L[u_p] |_{(p)} \big\|_{L^p(\Om)}\Big)^{1-p} \, {\av_\Om} \big| \L[u_p] \big|^{p-1}_{(p)}  \, \mathrm{d} \mL^n
\\
&\leq\, \Big(\big \| | \L[u_p] |_{(p)} \big\|_{L^p(\Om)}\Big)^{1-p} \left(\, \, {\av_\Om} \big|\L[u_p] \big|^{p}_{(p)}  \, \mathrm{d} \mL^n \right)^{\!\!\frac{p-1}{p}}
\\
& = \, 1.
\end{split}
\]
Next, fix $\phi \in C^2_c(\Om)$. Then, by using the regularity of $K$, we formally compute
\[
\begin{split}
\frac{\mathrm{d}}{\mathrm{d}\e}\bigg|_{\e=0} \!\! \E_p(u_p+\e\phi) &= 
p \! \left(\, \, {\av_\Ga} \big|\Q[u_p]-q^\de\big|^{p}_{(p)}  \, \mathrm{d} \mH^\ga \right)^{\!\!\frac{1}{p}-1} \!\! \, {\av_\Ga} \big|\Q[u_p]-q^\de\big|^{p-2}_{(p)}\big( \Q[u_p]-q^\de\big)   \centerdot
\\
& \hspace{78pt} \centerdot \Big[ K_r(\cdot,u_p,\D u_p) \,\phi \, +\, K_p(\cdot,u_p,\D u_p)\cdot \D \phi  \Big] \, \mathrm{d} \mH^\ga
\\
& \ \ \ +\, \al p \left(\, \, {\av_\Om} \big| \L[u_p] \big|^{p}_{(p)}  \, \mathrm{d} \mL^n\right)^{\!\!\frac{1}{p}-1} \!\! \, {\av_\Om} \big| \L[u_p] \big|^{p-2}_{(p)} \, \L[u_p] \, \L[\phi] \, \mathrm{d} \mL^n.
\end{split}
\]
Since $u_p$ is the minimiser of $\E_p$ in the space, we have that $\E_p(u_p)\leq \E_p(u_p+\e \phi)$ for all $\e\in\R$ and $\phi \in C^2_c(\Om)$. Therefore, this above computation implies that the PDE \eqref{1.19} is indeed satisfied as claimed in the statement of the lemma, upon confirming that the formal computation in the integrals above is rigorous, and that therefore $\E_p$ is Gateaux differentiable at the minimiser $u_p$ for any direction $\phi \in W^{2,p}_0(\Om)$. This is indeed the case: since $u_p \in (C^1 \cap W^{2,p})(\Om)$, $\L[u_p] \in L^p(\Om)$ and $\Q[u_p]-q^\de \in L^\infty(\Ga,\mH^\ga)$, the H\"older inequality implies that 
\[
\big| \L[u_p] \big|^{p-2}_{(p)} \, \L[u_p] \, \L[\phi] \ \in \, L^1(\Om)
\]
and 
\[
\big|\Q[u_p]-q^\de\big|^{p-2}_{(p)}\big( \Q[u_p]-q^\de\big) \Big[ K_r(\cdot,u_p,\D u_p) \,\phi + K_p(\cdot,u_p,\D u_p)\cdot \D \phi  \Big] \ \in \ L^1(\Ga,\mH^\ga),
\]
for any $\phi \in W^{2,p}_0(\Om) \sub C^1(\overline{\Om})$, because of the continuity of $K(x,r,p)$ in $x$ and the $C^1$ regularity in $(r,p)$.
\ep

\begin{lemma} \label{lemma4} For any $\al,\de>0$, consider the minimiser $u_\infty$ of $\E_\infty$ constructed in Lemma \ref{lemma2} as sequential limit of minimisers $(u_p)_{p>n}$ of the functionals $(\E_p)_{p>n}$ as $p_i\to \infty$. Then, there exist signed Radon measures $\mu_\infty \in \mM(\Om)$ and $\nu_\infty \in \mM(\Ga)$ such that the triplet $(u_\infty,\mu_\infty,\nu_\infty)$ satisfies the PDE \eqref{1.14} in the distributional sense, that is
\[
\int_\Om \Big( K_r(\cdot,u_\infty,\D u_\infty) \,\phi \, +\, K_p(\cdot,u_\infty,\D u_\infty)\cdot \D \phi \Big)\, \mathrm{d} \nu_\infty \, +\, \al \int_\Om \L[\phi] \, \mathrm{d} \mu_\infty\, =\, 0,
\]
for all $\phi \in C^2_c(\Om)$. Additionally, there exists a further subsequence along which the weak* modes of convergence of \eqref{1.18} hold true as $p\to \infty$.
\end{lemma}

\bp As noted in the beginning of the proof of Lemma \ref{lemma3}, we have the $p$-uniform total variation bounds $\|\mu_p\|(\Om)\leq 1$ and  $\|\nu_p\|(\Ga)\leq 1$. Hence, by the sequential weak* compactness of the spaces of Radon measures
\[
\mM(\Om)\, =  \big(C^0_0(\Om)\big)^*, \ \ \ \mM(\Om)\, =  \big(C^0(\Ga)\big)^*,
\]
there exists a further subsequence denoted again by $(p_i)_1^\infty$ such that $\mu_p \weakstar \mu_\infty$ in $\mM(\Om)$ and  $\nu_p \weakstar \nu_\infty$ in $\mM(\Ga)$, as $p_i \to \infty$. Fix now $\phi \in C^2_c(\Om)$. By Lemma \ref{lemma3}, we have that the triplet $(u_p,\mu_p,\nu_p)$ satisfies \eqref{1.19}, that is
\[
\begin{split}
\int_\Ga \Big( K_r(\cdot,u_p,\D u_p) \,\phi \, +\, K_p(\cdot,u_p,\D u_p)\cdot \D \phi \Big)  \, \mathrm{d} \nu_p \,+\, \int_\Om \L[\phi]  \, \mathrm{d} \mu_p\, =\, 0.
\end{split}
\]
Since 
\[
\L[\phi] \, \in \, C^0_0(\overline{\Om}), \ \ \  K_p(\cdot,u_\infty,\D u_\infty)\cdot \D \phi \, \in \, C^0(\Ga)
\]
and also
\[
\begin{split}
K_r(\cdot,u_p,\D u_p) \,\phi \, +\, & K_p(\cdot,u_p,\D u_p) \cdot \D \phi \larrow
\\
 &K_r(\cdot,u_\infty,\D u_\infty) \,\phi \, +\, K_p(\cdot,u_\infty,\D u_\infty)\cdot \D \phi,
\end{split}
\]
uniformly on $\Ga$ as $p_i \to \infty$ (as a consequence of the $C^1$ regularity of $K$ and the convergence $u_p\larrow u_\infty$ in $C^1(\overline{\Om})$), the weak*-strong continuity of the duality pairings between the above spaces of measures $\mM(\Om)$, $\mM(\Ga)$ and their respective predual spaces $C^0_0(\Om)$, $C^0(\Ga)$, allows us to conclude and obtain \eqref{1.14} by passing to the limit as  $p_i \to \infty$ in \eqref{1.19}.
\ep

\begin{remark} By testing in the weak formulation of \eqref{1.19} against $\phi \in C^2_c(\Om \set \Ga)$ (namely for those test functions such that $\phi \equiv 0$ on $\Ga$), we obtain $\L^*[\mu_\infty]= 0$ in $\Om\set \Ga$, that is
\[
\L^*\Big(\big| \L[u_p] \big|^{p-2}_{(p)} \, \L[u_p]\Big) \, =\, 0  \ \ \  \text{ in }\Om\set \Ga,
\]
in the distributional sense. Similarly, by testing in the weak formulation of \eqref{1.14} against $\phi \in C^2_c(\Om \set \Ga)$, we obtain
\[
\L^*[\mu_\infty] \, =\, 0  \ \ \  \text{ in }\Om\set \Ga,
\]
in the distributional sense.

\end{remark}

\begin{lemma} \label{lemma5} For any $\al,\de>0$, $p>n$ and $u^0 \in (W^{2,p}\cap W^{1,p}_g)(\Om)$ such that
\[
\big\| q^\de - \Q[u^0] \big\|_{L^\infty(\Ga,\mH^\ga)} \,\leq\, \de,
\]
the ($(\al,\de)$-dependent) minimiser $u_p$ of $\E_p$ (constructed in Lemmas \ref{lemma1}-\ref{lemma4}), satisfies the error bounds \eqref{1.24}, that is:
\[
\Big\| \Q[u_p] - \Q[u^0] \Big\|_{L^\infty(\Ga,\mH^\ga)}  \leq\, 2\de + \al \, \| \L[u^0]\|_{L^p(\Om)}.
\]
If additionally $u^0 \in \mathcal{W}^{2,\infty}_g(\Om)$, then the ($(\al,\de)$-dependent) minimiser $u_\infty$ of $\E_\infty$ (constructed in Lemmas \ref{lemma1}-\ref{lemma4}), satisfies the error bounds \eqref{1.20}, that is:
\[
\Big\| \Q[u_\infty] - \Q[u^0] \Big\|_{L^\infty(\Ga,\mH^\ga)}  \leq\, 2\de + \al \, \| \L[u^0]\|_{L^\infty(\Om)}.
\]
\end{lemma}

\bp Let us use the symbolisation $q^0:= \Q[u^0]$, noting also that $q^0 \in C^0(\Ga)$ and that we have the estimate 
\[
\|q^\de-q^0\|_{L^\infty(\Ga,\mH^\ga)} \, \leq \, \de.
\]
For any $p\in (n,\infty)$, the function $u_p$ is a global minimiser of $\E_p$ in $(W^{2,p}\cap W^{1,p}_g)(\Om)$. Therefore,
\[
\E_p(u_p) \, \leq\, \E_p(u^0).
\]
This implies the estimate
\[
\begin{split}
\big \|\Q[u_p]-q^\de \big\|_{L^p(\Ga,\mH^\ga)} & +\, \al \big \| \L[u_p] \big\|_{L^p( \Om)} 
\\
&\leq\, \big \|\Q[u^0]-q^\de \big\|_{L^p(\Ga,\mH^\ga)}+\, \al \big \| \L[u^0] \big\|_{L^p( \Om)}. 
\end{split}
\]
The latter estimate together with the Minkowski and H\"older inequalities, in turn yield
\[
\begin{split}
\big \|\Q[u_p]- \Q[u^0] \big\|_{L^p(\Ga,\mH^\ga)} & \leq\, \big \|\Q[u^0]-q^\de \big\|_{L^p(\Ga,\mH^\ga)}\,
\\
& \ \ \ \ +\, \big \|\Q[u^0]-q^\de \big\|_{L^p(\Ga,\mH^\ga)}+\, \al \big \| \L[u^0] \big\|_{L^p( \Om)}
\\
& = \, 2 \|q^\de-q^0\|_{L^\infty(\Ga,\mH^\ga)}  \, +\, \al \big \| \L[u^0] \big\|_{L^p(\Om)}
\\
& \leq 2\de + \al \, \| \L[u^0]\|_{L^p(\Om)},
\end{split}
\]
as claimed. To obtain the corresponding estimate for $u_\infty$ in the case that additionally $u^0 \in \mathcal{W}^{2,\infty}_g(\Om)$, we may pass to the limit as $p\to \infty$ in the last estimate above: indeed, consider the subsequence $p_i\to\infty$ along which we have the strong convergence $u_p\larrow u_\infty$ in $C^1(\overline{\Om})$ and therefore $\Q[u_p] \larrow \Q[u_\infty]$ uniformly on $\Ga$. Since by assumption $\L[u^0]\in L^\infty(\Om)$, the conclusion follows by letting $i\to \infty$ in the last estimate.
\ep

We now establish Proposition \ref{proposition11}.

\BPP \ref{proposition11}. {\bf (i)} Let $\mB^n_\rho(x)$ be the open $\rho$-ball of $\R^n$ centred at $x$. By the Lebesgue differentiation theorem (see e.g.\ \cite{FL}) applied to the measure $\nu\LL_X$ (namely to $\nu$ extended to $\R^n$ by zero on $\R^n \set X$) and by recalling that $\mB_\rho(x)$ symbolises the open ball in $X$, we have
\[
\begin{split}
f(x)\, & =\, \lim_{\rho \to 0} \bigg(\,{\av_{\mB^n_\rho(x)}}f\, \mathrm{d}(\nu\LL_X)\bigg)
\\
& =\, \lim_{\rho \to 0} \bigg(\frac{1}{\nu(\mB_\rho(x))}{\int_{\mB_\rho(x)}}f\, \mathrm{d}\nu\bigg)
\end{split}
\]
and therefore
\[
\begin{split}
f(x)\, & \leq\, \lim_{\rho \to 0} \bigg(\frac{1}{\nu(\mB_\rho(x))}{\int_{\mB_\rho(x)}}f\, \mathrm{d}\nu\bigg)
\\
&\leq\, \lim_{\rho \to 0} \bigg( \nu-\underset{\mB_\rho(x)}{\ess\,\sup}\, f \bigg)
\\
& =\, f^\bigstar(x),
\end{split}
\]
for $\nu$-a.e.\ $x\in X$.
\smallskip

\noi  {\bf (ii)} Fix $x\in X$ and $\e>0$. For any $\de \in (0,\e)$ and $y \in \mB_{\de}(x)$ we have the inclusion of balls
\[
\mB_{\e-\de}(y) \, \sub \, \mB_\e(x). 
\]
Hence, since the limit as $\e \to 0$ in the definition of $f^\bigstar$ is in fact an infimum over all $\e>0$, we have
\[
\begin{split}
\sup_{y \in \mB_\de(x)} f^\bigstar(y) \, &= \, \sup_{y \in \mB_\de(x)}\bigg[\lim_{\rho\to 0} \bigg(\nu-\underset{z \in \mB_\rho(y)}{\ess\,\sup} \, f(z)\bigg)\bigg]
\\
&= \, \sup_{y \in \mB_\de(x)}\bigg[\inf_{\rho> 0} \bigg(\nu-\underset{z \in \mB_\rho(y)}{\ess\,\sup} \, f(z)\bigg)\bigg]
\end{split}
\]
and therefore
\[
\begin{split}
\sup_{y \in \mB_\de(x)} f^\bigstar(y) \, & \leq \, \sup_{y \in \mB_\de(x)}\bigg[ \nu-\underset{z \in \mB_{\e-\de}(y)}{\ess\,\sup} \, f(z) \bigg]
\\
& \leq \, \sup_{y \in \mB_\de(x)}\bigg[ \nu-\underset{z \in \mB_{\e}(x)}{\ess\,\sup} \, f(z) \bigg]
\\
& = \,  \nu-\underset{z \in \mB_{\e}(x)}{\ess\,\sup} \, f(z) .
\end{split}
\]
By letting $\de\to 0$ and $\e\to 0$, we obtain
\[
\begin{split}
\lim_{\de\to 0} \bigg(\sup_{y \in \mB_\de(x)} f^\bigstar(y)\bigg) \,  \leq \, \lim_{\e\to 0} \bigg(\nu-\underset{z \in \mB_{\e}(x)}{\ess\,\sup} \, f(z)\bigg) \, =\, f^\bigstar(x) ,
\end{split}
\]
for any $x\in X$. Hence
\[
\limsup_{X \ni y\to x} f^\bigstar(y) \, \leq \, f^\bigstar(x), 
\]
for any $x\in X$, as desired.
\ms

\noi  {\bf (iii)} We begin by noting that for any $x\in X$ and $\e>0$ we have
\[
\nu-\underset{y \in \mB_\e(x)}{\ess\,\sup} \, f(y) \, \leq\, \nu-\underset{y \in X}{\ess\,\sup} \, f(y)
\]
which readily implies 
\[
\sup_{x\in X} f^\bigstar(x)\, =\, \sup_{x\in X} \bigg(\nu-\underset{y \in \mB_\e(x)}{\ess\,\sup} \, f(y)\bigg) \, \leq\, \nu-\underset{x \in X}{\ess\,\sup} \, f(x).
\]
Conversely, by the definition of the essential supremum, for any $\de>0$, the set
\[
X(\de) \,:=\, \bigg\{x\in X\ :\ f(x)  > \nu-\underset{y\in X}{\ess\,\sup} \, f(y) -\de \bigg\}
\]
satisfies
\[
\nu(X(\de))\,>\,0.
\]
By the Lebesgue-Besicovitch differentiation theorem (see e.g.\ \cite{FL}), $\nu$-a.e.\ point $x\in X_\de$ has density $1$, namely
\[
\lim_{\e \to 0} \frac{\nu\big(X(\de) \cap \mB^n_\e(x)\big)}{\nu(\mB^n_\e(x))}\, =\, 1, \]
where $\mB^n_\e(x)$ is the open $\e$-ball centred at $x$ with respect to $\R^n$. Hence, since
\[
\mB_\e(x) \, = \, X \cap \mB^n_\e(x),
\]
for any $\de>0$, there exists $x_\de \in X(\de)$ such that
\[
\nu\big(X(\de) \cap \mB_\e(x_\de)\big) \, =\, \nu\big(X(\de) \cap \mB^n_\e(x_\de)\big)\, >\, 0.
\]
Therefore, since
\[
\nu-\underset{y \in X}{\ess\,\sup} \, f(y) \, \leq\, \de \,+\, f(x), \ \ \ \nu-\text{a.e. }x\in X(\de),
\]
we deduce
\[
\begin{split}
\nu-\underset{y \in X}{\ess\,\sup} \, f(y) \, &\leq\, \de \,+\, \nu-\underset{y \in \mB_\e(x_\de) \cap X(\de)}{\ess\,\sup} \, f(y)
\\
&\leq\, \de \,+\, \nu-\underset{y \in \mB_\e(x_\de)}{\ess\,\sup} \, f(y).
\end{split}
\]
By letting $\e \to 0$ in the above inequality, we infer that
\[
\begin{split}
\nu-\underset{x \in X}{\ess\,\sup} \, f(x) \, & \leq\, \de \,+\, \lim_{\e\to 0}\bigg(\nu-\underset{y \in \mB_\e(x_\de)}{\ess\,\sup} \, f(y)\bigg)
\\
& = \, \de \,+\, f^\bigstar(x_\de)
\\
& \leq \, \de \,+\, \sup_{x\in X} f^\bigstar(x),
\end{split}
\]
for any $\de>0$. By letting $\de\to 0$, we obtain
\[
\nu-\underset{x \in X}{\ess\,\sup} \, f(x) \, \leq\, \sup_{x\in X} f^\bigstar(x),
\]
as desired. This inequality completes the proof.
\qed
\ms
\ms

By invoking Proposition \ref{proposition12} whose proof follows, we readily obtain \eqref{1.15}-\eqref{1.16} by choosing
\[
X=\Ga,\ \ \nu=\mH^\ga\LL_\Ga,\ \ f_k=\Q[u_{p_k}]-q^\de,\ \ f_\infty=\Q[u_{\infty}]-q^\de.
\]

\BPP \ref{proposition12}. {\bf (i)} By the definition of $\nu_k$, we have for any continuous function $\phi \in C^0(X)$ with $|\phi|\leq 1$ that
\[
\begin{split}
\left|\int_X\phi \, \mathrm{d} \nu_k \, \right| \, &\leq\, \frac{1}{ \big\| |f_k|_{(k)} \big\|^{k-1}_{L^k(X,\nu)}} \, {\av_X} \Big| \big(|f_k|_{(k)}\big)^{k-2} f_k \,\phi\Big|  \, \mathrm{d} \nu
\\
&\leq\, \frac{1}{ \big\| |f_k|_{(k)} \big\|^{k-1}_{L^k(X,\nu)}} \, {\av_X} \big(|f_k|_{(k)}\big)^{k-1}  \, \mathrm{d} \nu.
\end{split}
\]
Hence, by H\"older inequality, we have the total variation bound
\[
\begin{split}
\|\nu_k\|(X) \, &\leq\, \Big(\big \| | f_k |_{(k)} \big\|_{L^k(X,\nu)}\Big)^{1-k} \left(\, \, {\av_X} \big(|f_k|_{(k)}\big)^{k}  \, \mathrm{d} \nu \right)^{\!\!\frac{k-1}{k}}
\\
& = \, 1.
\end{split}
\]
By the sequential weak* compactness of the space $\mM(X)=\big(C^0(X) \big)^*$, we obtain the desired subsequence $(\nu_{k_i})_1^\infty \sub \mM(X)$ and the weak* sequential limit measure $\nu_\infty \in \mM(X)$.
\ms

\noi  {\bf (ii)} We begin by showing the elementary inequality
\[
\big| |f_k|_{(k)} - |f_\infty|\big| \, \leq \, \big| f_k - f_\infty\big| \,+\, \frac{1}{k}\ \ \text{ on }X.
\]
Indeed, if $|f_k|_{(k)} \geq |f_\infty|$, we have
\[
\begin{split}
\big| |f_k|_{(k)} - |f_\infty|\big| \, & =\, \sqrt{|f_k|^2 + k^{-2}}\, -\, |f_\infty|
\\
& \leq\, |f_k| - |f_\infty| \,+\, \frac{1}{k}
\\
& \leq \, \big| f_k - f_\infty\big| \,+\, \frac{1}{k}
\end{split}
\]
whilst if $|f_k|_{(k)} < |f_\infty|$, we have
\[
\begin{split}
\big| |f_k|_{(k)} - |f_\infty|\big| \, &=\, |f_\infty| \, -\, \sqrt{|f_k|^2 + k^{-2}}
\\
& \leq\,  |f_\infty| \,-\, |f_k|
\\ 
&\leq \, \big| f_k - f_\infty\big| \,+\, \frac{1}{k}.
\end{split}
\]
Fix now $\e>0$. The inequality we just proved implies that if $f_k \larrow f_\infty$ uniformly on $X$ as $k\to \infty$ (note that $f_k, f_\infty$ might be discontinuous), then $|f_k|_{(k)} \larrow |f_\infty|$ uniformly on $X$ as $k\to \infty$. Hence, there exists $k(\e)\in\N$ such that
\[
\|f_k -f_\infty\|_{L^\infty(X,\nu)}<\,\frac{\e}{4},\ \ \ \big\| |f_k|_{(k)} - |f_\infty|\big\|_{L^\infty(X,\nu)}<\,\frac{\e}{4},
\]
for all $k\geq k(\e)$. Therefore,
\[
\begin{split}
|f_k| \, &\leq\, |f_\infty| \, +\, \frac{\e}{4},\ \ \ \ \nu-\text{a.e. on }X,
\\
|f_k|_{(k)} \,& \geq\, |f_\infty| \, -\, \frac{\e}{4},\ \ \ \ \ \nu-\text{a.e. on }X.
\end{split}
\]
By integrating the latter inequality and using the Minkowski inequality, we obtain
\[
\| |f_k|_{(k)} \|_{L^k(X,\nu)} \, \geq \,\| f_\infty \|_{L^k(X,\nu)} \, -\, \frac{\e}{4},
\]
for all $k\geq k(\e)$. Since 
\[
\| f_\infty \|_{L^\infty(X,\nu)} \, =\, \lim_{k\to \infty} \| f_\infty \|_{L^k(X,\nu)} ,
\]
by choosing $k(\e)$ greater if necessary, we deduce
\[
\| |f_k|_{(k)} \|_{L^k(X,\nu)} \, \geq \,\| f_\infty \|_{L^\infty(X,\nu)} \, -\, \frac{\e}{2},
\]
for all $k\geq k(\e)$. Let now ${\mathrm{d} \nu_k}/{\mathrm{d}\nu}$ symbolise the Radon-Nikodym derivative of $\nu_k$ with respect to $\nu$. It follows that
\[
\frac{\mathrm{d} \nu_k}{\mathrm{d}\nu}\, =\,  \frac{1}{ \nu(X)} \frac{\big(|f_k|_{(k)}\big)^{k-2} f_k }{ \big\| |f_k|_{(k)} \big\|^{k-1}_{L^k(X,\nu)}}, \ \ \ \nu\text{-a.e.\ on }X.
\]
By the above, for any $\e>0$ small enough (recall that $f_\infty\not\equiv 0$) and for any $k\geq k(\e)$, we have the estimate
\[
\bigg| \frac{\mathrm{d} \nu_k}{\mathrm{d}\nu} \bigg|\, \leq\, \frac{1}{ \nu(X)} \left( \frac{ \ \ \dfrac{1}{k} +\, |f_\infty|\,+ \, \dfrac{\e}{4} \ \ }{ \| f_\infty \|_{L^\infty(X,\nu)} -\, \dfrac{\e}{2} } \right)^{\!\! k-1},\ \ \ \ \ \nu-\text{a.e. on }X.
\]
By choosing $k(\e)$ even larger if needed, we can arrange
\[
\bigg| \frac{\mathrm{d} \nu_k}{\mathrm{d}\nu} \bigg|\, \leq\, \frac{1}{ \nu(X)} \Bigg( \frac{ 2|f_\infty|\,+ \, \e  }{ 2\| f_\infty \|_{L^\infty(X,\nu)} -\, \e} \Bigg)^{\!\! k-1},\ \ \ \ \ \nu-\text{a.e. on }X.
\]
Since by Proposition \ref{proposition11} we have $|f_\infty| \leq |f_\infty|^\bigstar$ $\nu$-a.e.\ on $X$, we obtain
\[
\bigg| \frac{\mathrm{d} \nu_k}{\mathrm{d}\nu} \bigg|\, \leq\, \frac{1}{ \nu(X)} \Bigg( \frac{ 2|f_\infty|^\bigstar\,+ \, \e  }{ 2\| f_\infty \|_{L^\infty(X,\nu)} -\, \e} \Bigg)^{\!\! k-1},\ \ \ \ \ \nu-\text{a.e. on }X.
\]

Consider now for any $\e>0$ the $\nu$-measurable set
\[
X_\e\,:=\, \Big\{ |f_\infty|^\bigstar < \|f_\infty\|_{L^\infty(X,\nu)} -2\e\Big\}.
\]
Notice also that $X_\e$ is in fact open in $X$ because $|f_\infty|^\bigstar$ is upper semicontinuous (Proposition \ref{proposition11}). Additionally, we have the estimate
\[
\bigg| \frac{\mathrm{d} \nu_k}{\mathrm{d}\nu} \bigg|\, \leq\, \frac{1}{ \nu(X)} \Bigg( \frac{ 2\| f_\infty \|_{L^\infty(X,\nu)} - \, 3\e }{ 2\| f_\infty \|_{L^\infty(X,\nu)} -\, \e } \Bigg)^{\!\! k-1},\ \ \ \ \ \nu-\text{a.e. on } X_\e.
\]
The above estimate together with the Lebesgue Dominated Convergence theorem imply that for any $\e>0$ small enough we have
\[
\frac{\mathrm{d} \nu_k}{\mathrm{d}\nu} \larrow 0\ \ \text{ in }L^1(X_\e,\nu),\ \text{ as }k\to\infty.
\]

Consider now the sequence of nonnegative total variation measures $(\|\nu_k\|)_1^\infty \sub \mM(X)$. Since this sequence is also bounded in the space, there exists a nonnegative limit measure $\la_\infty$ such that
\[
\|\nu_k\| \, \weakstar \, \la_\infty \ \ \text{ in } \mM(X),
\]
along perhaps a further subsequence $(k_i)_1^\infty$. Additionally, since $\nu_k \weakstar \nu_\infty$ in $\mM(X)$, we have the inequality (see e.g. \cite{AFP})
\[
\| \nu_\infty \| \, \leq\, \la_\infty.
\]
Note now that for each $k\in\N$, by the Lebesgue-Radon-Nikodym theorem applied to $\|\nu_k\|<<\nu$ we have the decomposition
\[
\|\nu_k\| \, = \bigg|\frac{\mathrm{d} \nu_k}{\mathrm{d}\nu}\bigg|\, \nu.
\]
Hence, we infer that
\[
\|\nu_k\|(X_\e)\, \leq \, \int_{X_\e} \bigg|\frac{\mathrm{d} \nu_k}{\mathrm{d}\nu} \bigg|\,\mathrm{d}\nu\, \larrow\, 0, \ \ \ \text{ as }k\to\infty.
\]
Therefore, since $X_\e$ is open in $X$, by the weak* lower-semicontinuity of measures on open sets (see e.g.\ \cite{FL,AFP}) and the above arguments, we have
\[
\begin{split}
\|\nu_\infty\|(X_e)\, & \leq\, \la_\infty(X_\e) 
\\
& \leq\, \liminf_{i\to \infty} \|\nu_{k_i}\| (X_\e)
\\
& \leq\, \liminf_{i\to \infty} \int_{X_\e} \bigg|\frac{\mathrm{d} \nu_{k_i}}{\mathrm{d}\nu} \bigg|\,\mathrm{d}\nu
\\
& =\, 0.
\end{split}
\]
Therefore, we have obtained
\[
\nu_\infty \Big(\Big\{ |f_\infty|^\bigstar< \|f_\infty\|_{L^\infty(X,\nu)}-2\e \Big\}\Big)\, =\, 0, \ \ \text{ for any }\e>0.
\]
By letting $\e\to 0$ along the sequence $\e_j:=2^{-j-1}$, the continuity of the measure $\nu_\infty$ implies
\[
\begin{split}
\nu_\infty \Big( \Big\{ |f_\infty|^\bigstar < \|f_\infty\|_{L^\infty(X,\nu)} \Big\} \Big)\, &=\, \nu_\infty \Bigg( \bigcup_{j=1}^\infty\Big\{ |f_\infty|^\bigstar< \|f_\infty\|_{L^\infty(X,\nu)}-2^{-j} \Big\}\Bigg)
\\
&=\, \lim_{j\to \infty}\nu_\infty \Big(\Big\{ |f_\infty|^\bigstar< \|f_\infty\|_{L^\infty(X,\nu)}-2^{-j} \Big\}\Big)
\\
&=\, 0.
\end{split} 
\]
Then, the definition of support of the measure $\nu_\infty$ and the upper semicontinuity of the function $|f_\infty|^\bigstar$ on $X$ (by Proposition \ref{proposition11}) yield
\[
\begin{split}
X \set \mathrm{supp}(\nu_\infty) \, &= \, \bigcup \Big\{U \sub X\text{ open }:\ \nu_\infty(U)=0 \Big\}
\\
& \supseteq \, \Big\{ |f_\infty|^\bigstar< \|f_\infty\|_{L^\infty(X,\nu)} \Big\}.
\end{split}
\]
In conclusion, we infer that
\[
\begin{split}
\mathrm{supp}(\nu_\infty) \, & \sub \, X\set \Big\{ |f_\infty|^\bigstar< \|f_\infty\|_{L^\infty(X,\nu)} \Big\}  \phantom{\bigg|}
\\
& =\, \Big\{ |f_\infty|^\bigstar= \|f_\infty\|_{L^\infty(X,\nu)} \Big\},
\end{split}
\]
as desired. 

\ms

\noi {\bf (iii)} Suppose that $|f_\infty|$ is continuous on $X$ and recall the properties of the essential limsup established in Proposition \ref{proposition11}. Then, for any $x\in X$ we have 
\[
\begin{split}
\Big| |f_\infty|^\bigstar(x) - \, |f_\infty|(x)\Big|\, &=\, \left|\lim_{\e\to 0}\bigg(\nu-\underset{\mB_\e(x)}{\ess\sup}\, |f_\infty|\bigg)  - \, |f^\infty|(x) \right|
\\
&=\, \left|\lim_{\e\to 0}\bigg(\nu-\underset{\mB_\e(x)}{\ess\sup}\, |f_\infty| \, - \, |f^\infty|(x)\bigg) \right|
\\
&\leq \, \limsup_{\e \to 0}\left|  \nu-\underset{\mB_\e(x)}{\ess\sup}\, \Big(|f_\infty|\,  - \, |f^\infty|(x)\Big) \right|
\\
 &\leq \, \limsup_{\e \to 0}\Big\||f_\infty| \, - \, |f_\infty|(x)\Big\|_{L^\infty(\mB_\e(x),\nu)}
\\
&=\, 0,
\end{split}
\]
showing that $|f_\infty|^\bigstar \equiv |f_\infty|$, if it holds that $|f_\infty|$ is continuous on $X$.
\qed

\ms

\ms

\noi \textbf{Acknowledgement.} The author would like to thank Jochen Broecker for discussions on inverse source identification problems, as well as Roger Moser, Jan Kristensen and Tristan Pryer for inspiring scientific discussions on the topics of Calculus of Variations in $L^\infty$. He is also indebted to the anonymous referee for their constructive comments which improved the content and the presentation of an earlier version of this paper.

\ms

\bibliographystyle{amsplain}

\end{document}